\documentclass[amssymb,amscd,12pt]{article}
\usepackage{epsfig}
\usepackage{amsmath}
\usepackage{amsthm}
\usepackage{amssymb}

\usepackage{amscd}
\usepackage{amsfonts}  
\usepackage{makeidx}
\usepackage{hyperref}
\usepackage[all]{xy}
\usepackage{graphicx}    
\usepackage{multicol}    

\newtheorem{theorem}{Theorem}

\newtheorem{proposition}[theorem]{Proposition}

\newtheorem{remark}[theorem]{Remark}
\newtheorem{corollary}[theorem]{Corollary}

{\catcode`@=11\global\let\c@equation=\c@theorem}

\newcommand{\id}{\operatorname{id}}

\title {On the classification of $1$-connected $7$-manifolds with torsion free second homology}
\author{Matthias Kreck }
\date{}

\begin{document}
\maketitle

{\small{\centerline {Abstract} We generalize a result (\cite{Kreck}, Theorem 6) of the author about the classification of $1$-connected $7$-manifolds and demonstrate its use by two concrete applications,  one to $2$-connected $7$-manifolds  (a new proof -- and slightly different formulation -- of an up to now unpublished Theorem by Crowley and Nordstr\"om \cite{C-N}) and one to simply connected $7$-manifolds with the cohomology ring of $S^2 \times S^5 \sharp S^3 \times S^4$. The answer is in terms of generalized Kreck-Stolz invariants, which in the case of $2$-connected $7$-manifolds is equivalent to a quadratic refinement of the linking form and a generalized Eells-Kuiper invariant.}}

\section{Introduction} We generalize a result (\cite{Kreck}, Theorem 6) about the classification of $1$-connected $7$-manifolds  and demonstrate its use by two concrete applications. The first application of our main theorem (Theorem 5) is a new proof (and slightly different formulation) of an up to now unpublished Theorem by Crowley and Nordstr\"om \cite{C-N} about the classification of $2$-connected $7$-manifolds. To formulate the theorem in a convenient way we define the concept of a {\bf d-structure}. Let $M$ be a closed $2$-connected  $7$-manifold with spin-structure. Then it has a spin Pontrjagin class $\bar p_1(M)$ (the pullback of the universal spin Pontrjagin class $\bar p_1:= p_1/2 \in H^4(BSpin)$, where $p_1 \in H^4(BSpin)$ is the pull back of $p_1 \in H^4(BSO)$ under the map induced by the canonical projection $BSpin \to BSO$). The divisibility of $\bar p_1(M)$ in $H^4(M)/_{\text {torsion}}$ is denoted by $d(M)$ (if $\bar p_1(M)$ is torsion, we set $d(M) =0$). For a spin $7$-manifold it was proven in \cite{K-S2}, Lemma 6.5, that the spin Pontrjagin class $\bar p_1(M)$ reduces mod $2$ to the $4$-th Stiefel-Whitney class, which for spin-manifolds is the $4$-th  Wu-class and by definition is $0$ for dimensional reasons (see \cite{C-N}, Lemma 2.2 (i)). Thus $d(M)$ is even. 

A {\bf $d(M)$-structure} is a class $x \in H^4(M)$ with $\bar p_1(M) - d(M) x$ torsion, where we take $x= 0$ if $d(M)=0$. A $d(M)$-structure can be used to define two invariants of $M$, previously defined by Crowley and Nordstr\"om, naturally occurring as Kreck-Stolz type invariants, when we apply our main theorem (Theorem 5): a {\bf quadratic refinement} $q_x: \text{torsion}\,\, H^4(M) \to \mathbb Q/_{2\mathbb Z}$ of the linking form $b: \text {torsion}\,\, H^4(M)\otimes \text{torsion}\,\, H^4(M) \to \mathbb Q/_{\mathbb Z}$ and an invariant called the {\bf generalized Eells-Kuiper}  invariant  $GEK(M,x)\in \mathbb Q/_{8 \cdot gcd( 28, \frac {d(M)} 2,\frac{ d(M)^2+2d(M)}8)\mathbb Z}$. We will give the details of these invariants in section 3. 

Thus we obtain the following data:
$$
\varphi (M,x) := (H^4(M), \bar p_1(M) , q_x(M), GEK(M,x)).
$$
We say that two quadruples $\varphi(M,x)$ and $\varphi(M',x')$   are isomorphic, if there is an isomorphism between the cohomology groups mapping $x$ to $x'$, mapping  $\bar p_1(M)$ to $\bar p_1(M')$, commuting with the quadratic refinements and such that $GEK(M,x) = GEK(M',x')$. The group $G$ of automorphisms of $H^4(M)$ preserving $\bar p_1(M)$ acts on the quadruples and the equivalence class is called $\varphi(M)$.

\begin{theorem} (see \cite{C-N}, Theorem 1.3, see also Theorem 1.4) Let $M$ and $M'$ be  $2$-connected closed oriented $7$-manifolds with $d(M) = d(M') =d$. Let $x\in H^4(M)$ and $x' \in H^4(M')$ be $d$-structures. Then $M$ and $M'$ are $d$-structure and  orientation preserving diffeomorphic if and only if $\varphi(M,x)$ is isomorphic to $\varphi(M',x')$. There is a diffeomorphism inducing a prescribed isomorphism.

As a consequence $M$ is orientation preserving diffeomorphic to $M'$ if and only $\varphi(M) = \varphi(M')$. 

\end{theorem}

We will later give formulas for the change of the invariant $\varphi(M,x)$ if we change the $d(M)$-structure (Proposition 7). In \cite{C-N}, it was shown that the invariants used there can be realized. 

The second application concerns closed simply connected spin 
$7$-manifolds $M$ with the cohomology ring of $S^2 \times S^5 \sharp S^3 \times S^4$. Again we consider $d(M)$, the divisibility of $\bar p_1(M)$ and recall that $d(M) $ is even. Using the characteristic numbers occurring in our main theorem we define Kreck-Stolz type invariants 
$$
s_1(M) \in \mathbb Z/_{8 \cdot gcd( 28, \frac {d(M)} 2,\frac{ d(M)^2+2d(M)}8)\mathbb Z}
$$
$$
s_2(M) \in \mathbb Z/_{gcd(24,d(M))\mathbb Z}
$$
and
$$
s_3(M) \in \mathbb Z/_{2\mathbb Z}.
$$
We will give the definitions in section 4. 

\begin{theorem} Two closed oriented simply connected spin $7$-manifolds $M$ and $M'$ with same cohomology ring as $S^2 \times S^5 \sharp S^3 \times S^4$ are orientation preserving diffeomorphic  if and only if $d(M) = d(M')$ and
$$
s_i(M) = s_i(M')
$$
for $i = 1,2,3$.
\end{theorem}

It is easy to find, for a given integer $s$, a manifold $M$ with the cohomology ring of 
$S^2 \times S^5 \sharp S^3 \times S^4$  such that $d(M) = 2s y$ for some generator $y \in H^4(M)$.  For this we consider an $S^3$-bundle over $S^4$ with Euler class $0$ and first Pontrjagin class $\pm 4s y$ (recall that $\pi_3(SO(4) \cong \mathbb Z \oplus \mathbb Z$ and that the isomorphism is given by the Euler class $e(E)$ and $(p_1(E)-2e(ME))/4$, where $E$ is the corresponding vector bundle) and take the connected sum of the total space, denoted by $S_s$, with $S^2 \times S^5$. This has an orientation reversing self diffeomorphism (since the bundle has a section). More generally, consider the total space of the sphere bundle over $\mathbb CP^2$ of a $4$-dimensional vector bundle with second Stiefel-Whitney class non-zero, Euler class $x^2$ and first Pontrjagin class $(3 + 4k)x^2$, where $x \in H^2(\mathbb CP^2)$ is a generator. We orient the vector bundle such that the self intersection number on the middle homology is $+1$ and take the induced orientation on the boundary. Such bundles can be obtained by starting from the sphere bundle of the Whitney sum of  the tautological line bundle over $\mathbb CP^2$ with the trivial complex line bundle and reattaching the bundle over the top cell appropriately. The total space is denoted by $T_k$ and it has the cohomology ring of $S^2 \times S^5$ (this follows from the Gysin sequence). Thus the connected sum
$$
M_{s,k} := S_s \sharp T_k
$$
is a manifold of the type we consider. Furthermore, if $\Sigma_r$ is a $7$-dimensional exotic sphere, which is given by the $r$-fold connected sum of the boundary of the $E_8$-plumbing in dimension $8$,  also
$$
M_{s,k} \sharp \Sigma_r
$$
is in this class. When are two such manifolds diffeomorphic? This can easily be decided using Theorem 2 by computing the invariants $s_i(M_{s,k}\sharp \Sigma_r)$. We will carry this out in section 5 and prove the following corollary.

\begin{corollary} The manifolds $M_{s,k} \sharp \Sigma_r$ and $M_{s',k'} \sharp \Sigma_{r'}$ are diffeomorphic if and only if 
$$
s = s'
$$
$$
4(1-r)  + 3k + 2k^2 = 4(1-r')  + 3k' + 2(k')^2 \,\, {\text mod} \,\,4\,\,  \gcd (28,\frac{s(s-1)}2,s^2)
$$
and 
$$
k = k' \,\, \text {mod} \,\,gcd( 12,s)
$$

\end{corollary}

The formula allows the determination of the inertia group of $M_{s,k}$. The inertia group of a smooth $n$-manifold $M$ consists of the subgroup of the exotic spheres $\Sigma$ such that $\Sigma \sharp M$ is diffeomorphic to $M$. In dimension $7$ the group of exotic spheres is isomorphic to $\mathbb Z/_{28}$ and the isomorphism is given by the Eells-Kuiper invariant, which for exotic spheres agrees with the generalized Eells-Kuiper invariant. The formula in the corollary implies that $\Sigma_r \sharp M_{s,k} $ is diffeomorphic to $M_{s,k}$ if and only if $r = 0 \,\, \text{mod}\,\,\gcd (28,\frac{s(s-1)}2,s^2)$. Thus all subgroups of $\Theta_7$ occur as inertia groups of appropriate $7$-manifolds with cohomology ring of $S^2 \times S^5 \sharp S^3 \times S^4$.\\

Now we explain our main theorem. Since the formulation for non-spin manifolds is a bit technical we restrict ourselves here to the spin case. For integers $n$ and $m$ we denote by $B(n,m)$ the fibration over $BSO$ with total space $K(\mathbb Z^n,2) \times K(\mathbb Z^m,4) \times BSpin$ and map given by the projection to $BSpin$ composed with the projection $BSpin \to BSO$. If $M$ is a closed 1-connected $7$-dimensional spin manifold with torsion free second homology we consider a lift $\bar \nu: M \to B(n,m)$ of the normal Gauss map and require that  the induced map  $H^2(B(n,m)) \to H^2(M)$ is an isomorphism and the induced map $H^4(B(n,m)) \to H^4(M)$ is surjective. Note that $\bar \nu$ is equivalent to a spin structure on $M$ (which is unique, since $M$ is 1-connected), an isomorphism $\alpha: \mathbb Z^n \to H^2(M)$ and a homomorphism $\beta : \mathbb Z^m \to H^4(M)$, such that the image of $\beta$ together with the spin Pontrjagin class $\bar p_1(M)$  and products of elements in $H^2(M)$ generate $H^4(M)$. Since the spin-structure on $M$ is unique, we can omit it from the data. We call the pair $(\alpha, \beta)$ (or the corresponding lift $\bar \nu(\alpha,\beta)$) with these properties a {\bf polarization} of $M$ in $B(n,m)$. We note here, that $B(n,0)$ is the normal $2$-type of $M$, which we enrich by the product with $K(\mathbb Z^m,4)$ and $\bar \nu(\alpha,\beta)$ is an enriched normal $2$-smoothing. We will discuss and generalize this in section 2. 

Let $(M', \bar \nu(\alpha',\beta'))$ be another closed spin manifold with torsion free second homology together with a polarization, a pair of isomorphisms $g: H^2(M') \to H^2(M)$ and $h: H^4(M') \to H^4(M)$ is called a {\bf multiplicative tangential isomorphism}, if $h$ preserves the first spin Pontrjagin class (which in our situation determines the stable tangent bundle)  and  $g(x) \cup g(y) = h(x \cup y)$ for all $x,y \in H^2(M)$, and the map $g$ commutes with $\alpha$ and  $\alpha'$, and $h$ commutes with   $\beta$ and $\beta'$. The existence of a polarizations and a multiplicative tangential isomorphism is a necessary condition for the existence of a diffeomorphism from $M$ to $M'$ compatible with the polarizations.

\begin{theorem} Let $M$  and $M'$  be simply connected closed spin $7$-manifolds with torsion free second homology. 

Then if $M'$ is another closed  simply connected $7$-manifold, $M'$  is orientation preserving diffeomorphic to $M$ if  and only if there  are polarizations $(\alpha, \beta)$ and $(\alpha', \beta')$ respectively and a multiplicative tangential isomorphism $(g,h)$ and a $B$ bordism $(W,l)$ between $(M, \bar \nu(\alpha, \beta)$ and $(M', \bar \nu(\alpha', \beta')$ such that\\
- $\text {sign} (W) =0$\\
- $\langle l^*(x) \cup l^*(y), [W,\partial W]\rangle  = 0$ for all classes $x$ and $y$ in $H^4(B; \mathbb Q)$ that map to zero in $H^4(\partial W; \mathbb Q)$. \\
Moreover there is a diffeomorphism $f: M \to M'$ that induces $(g,h)$ in (co)homology. 

The same statements hold if we replace the manifolds and bordisms by topological manifolds and diffeomorphism by homeomorphism.

\end{theorem}

Here the expression $\langle l^*(x) \cup l^*(y), [W,\partial W]\rangle  = 0$ has to be understood by pulling the classes back to $H^4(W,\partial W; \mathbb Q)$, taking the cup product and evaluating on the fundamental class. The relation between $(W,l)$ and $(g,h)$ is the following. Let $i_M$ and $i_{M'}$ be the inclusions, then $l^* i_{M}^* g = l^* i_{M'}^*: H^2(M') \to H^2(B)$ and $l^* i_{M}^* h = l^* i_{M'}^*: H^4(M') \to H^4(B)$. 

If one wants to apply the theorem to the classification of manifolds the first obstruction is the existence of a $B(n,m)$ bordism between $(M,\bar \nu(\alpha,\beta))$ and $(M',\bar \nu(\alpha',\beta'))$. We will show that the $7$-dimensional bordism group is zero, so that such a bordism always exists. Then one has to choose one and apply Theorem 4. If the invariants for the chosen bordism $(W,l)$ don't vanish, that doesn't mean that $M$ and $M'$ are not diffeomorphic. Then one has to alter the bordism, which up to bordism corresponds to the disjoint union with a closed $8$-dimensional $B$-manifold. Thus one has also to compute the $8$-dimensional bodsim group, which we will do for our applications. \\

The theorem is a generalization  of Theorem 6 in \cite{Kreck}, which is the special case, that $H^4(M)$ is finite and $m=0$, i.e. the $4$-th cohomology is generated by products of $2$-dimensional classes and $\bar p_1(M)$. This result had numerous applications, in particular to the classification of homogeneous spaces. Such applications are the motivation for the new result, since some simply connected $7$-manifolds have interesting geometric structures like metrics with positive curvature, Einstein metrics or $G_2$-structures. For those manifolds (and in general) explicit classification results in terms of numerical invariants are important. 

We will also prove a very special case of the main theorem for manifolds with boundary (see Theorem 5). This can be applied to study the mapping class group of simply connected $6$-manifolds with torsion free second homology group. Almost nothing was known about this. In joint work with Su Yang \cite{K-Y} we apply the theorem to determine the mapping class group for a large class of $6$-manifolds including examples of Calabi-Yau 3-folds.

The author would like to thank Diarmuid Crowley for stimulating discussions and Stephan Stolz, Peter Teichner and Su Yang for useful comments. Particular thanks to the referee for an unusual careful reading and numerous suggestions (including pointing at an error in the formulation of Theorem 4) which helped to improve the paper considerably.

\section{Main theorem} 

We assume that the reader is familiar with the basic definitions and results of modified surgery \cite{Kreck}. The definitions needed are: Normal $k$-type, normal smoothings, definition of $l$-monoid, and elementary obstructions. The main result needed is that if the obstruction in the $l$-monoid is elementary, surgery to an $s$-cobordism is possible \cite{Kreck}, Theorem 3.

We want to generalize Theorem 4 to the non-spin case and to the case where $M$ is of the form $N \times [0,1]$ for some 1-connected closed $6$-manifold $N$.

We begin with some notation and construction. For a natural number $n$  and a class $w_2\in H^2(K(\mathbb Z^n,2);\mathbb Z/2)$ we denote by $B(n,w_2)$ the fibration over $BSO$ which is the pullback of the fibration $K(\mathbb Z^n,2)  \to K(\mathbb Z/2,2)$ given by $w_2$ under the map $BSO \to K(\mathbb Z/2,2)$ given by the universal second Stiefel-Whitney class. 

\[
\xymatrix@=30pt{
B(n,w_2) \ar@{^-->}[d]^{p}  \ar@{^-->}[r]^{} & K(\mathbb Z^n,2)  \ar@{^-->}[d]^{w_2}  \\
BSO \ar@{^-->}[r]^{w_2} & K(\mathbb Z/2,2) \\}
\]

Let $M$ be a $1$-connected $7$-manifold with torsion free second homology. We consider an isomorphism $\alpha: \mathbb Z^n \to H^2(M) $. By abuse of notation we denote the homotopy class of maps $\alpha: M \to K(\mathbb Z^n,2)$ inducing $\alpha$ in cohomology with the same name. We denote the image of $w_2(M)$ under $(\alpha ^*) ^{-1}$ by $w_2(\alpha)$. The normal $2$-type of $M$ is $B(n,w_2(\alpha))$. Since $M$ is simply connected the map $\alpha$ determines a normal $2$-smoothing $\bar \nu_\alpha$. Namely, $\alpha$ gives a map $M \to K(\mathbb Z^n,2)$ inducing an isomorphism on second homology and the normal Gauss map is  a map to $BSO$, such that the composition of the two maps with the  maps to $K(\mathbb Z/2,2)$ given by $w_2(\alpha)$ and $w_2 \in H^2(BSO;\mathbb Z/2)$ are homotopic. By the definition of a pull back square we obtain a normal $2$-smoothing $\bar \nu_\alpha$. 

As we have seen in the spin case in the introduction we enrich the normal $2$-type and normal $2$-smoothing by additional data. This is a map $\beta: \mathbb Z^m \to H^4(M) $. Again we use the same letter for the corresponding map $\beta : M \to K(\mathbb Z^m,4)$. We denote the product of $B(n,w_2)$  with $K(\mathbb Z^m,4)$ by $B(n,m,w_2) := B (n,w_2) \times K(\mathbb Z^m,4) \to BSO$, where the map is the composition of the projection to $B(n,w_2)$ with the fibration above. The normal $2$-smoothing together with $\beta$ determines an enriched normal $2$-smoothing $\bar \nu(\alpha, \beta)$ with values in $B(n,m,w_2(\alpha))$. 

We call $(\alpha, \beta)$ a {\bf polarization of $M$} if $\alpha$ is an isomorphism and the induced map $H^4(B(n,m,w_2(\alpha)) \to H^4(M)$ is surjective.

Let $(M', \bar \nu(\alpha',\beta'))$ be another  $7$-manifold with torsion free second homology together with a polarization. A pair of isomorphisms $g: H^2(M') \to H^2(M)$ preserving the second Stiefel-Whitney class considered as homomorphism to $\mathbb Z/2$, and $h: H^4(M') \to H^4(M)$ is called a {\bf multiplicative tangential isomorphism}, if   $g(x) \cup g(y) = h(x \cup y)$ for all $x,y \in H^2(M')$, and the map $g$ commutes with $\alpha$ and  $\alpha'$, and $h$ commutes with   $\bar \nu(\alpha, \beta)^*$ and $\bar \nu(\alpha', \beta')^*$. The existence of a multiplicative tangential isomorphism is a necessary condition for the existence of a diffeomorphism from $M$ to $M'$ compatible with the polarizations. If $(W,l)$ is a bordism between $(M, \bar \nu(\alpha,\beta))$ and $(M', \bar \nu(\alpha',\beta'))$ we look for necessary and sufficient conditions so that the surgery obstruction in $l_8(1)$ occurring in \cite{Kreck}, Theorem 3, is elementary, which implies that $(W,l)$ is bordant rel. boundary to an $h$-cobordism. 

So far, we have not assumed that the manifolds are closed. If $M$ and $M'$ are compact manifolds with boundary and $h: \partial M \to \partial M'$ is a diffeomorphism we ask for an extension $f$ of $h$. Although we expect that our methods can be applied to more general manifolds with non-empty boundary we  only study here a very special but particularly interesting case, namely where  $M$ is of the form $N \times I$ with $N$ a simply connected $6$-manifold with $H_2(M)$ torsion free. The reason for looking at this case is that if $h: N \to N$ is an orientation preserving self diffeomorphism, we want to decide whether $h$ can be extended to a diffeomorphism on $N \times I$, which on the other boundary component restricts to the identity map or equivalently whether $h$ is pseudo-isotopic to the identity. By Cerf \cite{Ce}, Theorem 0, this implies that $g$ is isotopic to the identity. As above for $M$ we consider a polarization $(\alpha, \beta)$ of $N$ or equivalently of $N  \times I$ and the corresponding normal $2$-smoothing $\bar \nu(\alpha, \beta)$. A necessary condition for $h$ to be isotopic to the identity is that $h$ induces the identity on all homology groups. Another necessary condition is the following. If $h$ is isotopic to the identity, then the mapping torus $N_h$ is, as a fibre bundle over $S^1$, diffeomorphic to $S^1 \times N$, and we can choose a normal structure $\bar \nu_{N_h}$ on the mapping torus, whose restriction to the fibre $N$ is $\bar \nu(\alpha, \beta)$, the image of $(\bar \nu_{N_h})_* : H_4(N_h) \to H_4(B)$ is equal to the image of $(\bar \nu_{N_h}|_N)_*: H_4(N) \to H_4(B)$, and $(N_h,\bar \nu_{N_h})$ is zero bordant in $B$. We call such a normal $B$-structure on the mapping torus $N_h$  a {\bf good normal $B$-structure extending $\bar \nu(\alpha, \beta)$}. If $(W,l)$ is a zero bordism of $(N_h,\bar \nu_{N_h})$ we can ask whether it is $B$-bordant to a relative $h$-cobordism, between $N \times I$ and $N \times I$. If this is the case, the relative $h$-cobordism theorem implies that $h$ is pseudo isotopic and so isotopic to the identity. 

Given $(W,l)$ as in both cases above we can consider the following generalized relative characteristic numbers. Let $x$ and $y$ be in $H^4(B;\mathbb Q)$ such that the restriction of $l^*(x)$ and $l^*(y)$ vanish in $H^4(\partial W;\mathbb Q)$. Then
$$
\langle  l^*(x) \cup l^*(y), [W,\partial W]\rangle  \in \mathbb Q, 
$$ 
which is defined and well defined since $l^*(x)$ and $l^*(y)$ are in the image of $H^4(W,\partial W; \mathbb Q) \to H^4(W;\mathbb Q)$ and so the rational numbers $\langle  l^*(x) \cup l^*(y), [W,\partial W]\rangle  \in H^4(W, \partial W;\mathbb Q)$ can be computed. It is well defined since if we add $d(z)$ for some $z \in H^3\partial W$ to $x$ (or similarly to $y$), where $d$ is the boundary operator in the long pair sequence, the $d(z) \cup y = d(z \cup i_*(y)) =0$, where $i$ is the inclusion $\partial W \to W$.

The following Theorem is a generalization of Theorem 4 including the non-spin case. 

\begin{theorem} 
~\\
a) Let $M$ be a simply connected compact $7$-manifold with torsion free second homology.

Then if $M'$ is another closed  simply connected $7$-manifold, $M'$  is orientation preserving diffeomorphic to $M$ if and only there are polarizations $(\alpha, \beta)$ and  $(\alpha', \beta')$  respectively and a multiplicative tangential isomorphism $(g,h)$ and a $B$ bordism $(W,l)$ between $(M, \bar \nu(\alpha, \beta))$ and $(M', \bar \nu(\alpha', \beta'))$ such that \\
- $\text {sign} (W) =0$\\
- $\langle l^*(x) \cup l^*(y), [W,\partial W]\rangle  = 0$ for all classes $x$ and $y$ in $H^4(B; \mathbb Q)$ which map to zero in $H^4(\partial W; \mathbb Q)$. \\
Moreover there is a diffeomorphism $f: M \to M'$ which induces $(g,h)$ in (co)homology. 

b) Let  $h$ is an orientation preserving self-diffeomorphism of a closed simply connected $6$-manifold $N$ with torsion free second homology. We choose a polarization $(\alpha, \beta)$ of $N$. Then $h$ is isotopic to the identity if and only if $h$ acts by identity on all homology groups and there is a good normal structure $\bar \nu_{N_h}$ extending $\bar \nu(\alpha, \beta)$ and a zero bordism  $(W,l)$ over $B$ of it,  such that the same two conditions hold as in the closed case. 

The same statements hold if we replace the manifolds and bordisms by topological manifolds and diffeomorphism by homeomorphism.

\end{theorem}

To apply this theorem one has the freedom to choose the polarization resp. the good normal  structure. The only condition which is needed in the proof of Theorem 5  is a polarization $(\alpha, \beta)$, which in the case of $\beta$ gives a freedom in the choice of $m$. One should choose these structures with minimal number $m$. Otherwise the invariants are overdetermined. 

One might wonder why it is not required that the multiplicative tangential isomorphism respects the linking forms $b$  on the torsion of $H^4(M)$ and $H^	4(M')$. The reason is that this is automatically controlled. Namely the linking form on the torsion of $H^4(M)$  is determined by characteristic numbers, since if $x ,y \in H^4(M)$ are torsion elements, then there are $\hat x , \hat y \in H^4(B)$ mapping to $x$ and $y$ respectively. Then $b(x,y) = \langle l^*(\hat x) \cup l^*(\hat y) ,[W,\partial W]\rangle  \in \mathbb Q/_{\mathbb Z}$. This implies that under our conditions the linking form on the torsion subgroup of $H^4(M)$ and $H^4(M')$ is automatically preserved under a multiplicative tangential isomorphism.

This theorem is a generalization of Theorem 6 in \cite{Kreck}. There we considered a special class of manifolds with torsion free second homology and $4$-th cohomology torsion generated by products of classes in the second cohomology and the (twisted) $spin$-Pontrjagin class. Thus in this case the control given by $ \alpha$ is enough since we have a surjection $H^4(B) \to H^4(M)$. 

\begin{proof} The proof for topological and smooth manifolds is identical, except that at the end we apply the smooth or topological $h$-cobordism theorem. The necessity of the conditions is clear since, if $W$ is an $h$-cobordism, $\text {sign}(W) =0$ and the cup product condition holds automatically. Thus we assume that the conditions are fulfilled and prove that then there is a relative $h$-cobordism such that the $h$-cobordism theorem gives the desired diffeomorphism. 

 We apply the modified surgery theory of \cite{Kreck}. By \cite{Kreck}, Proposition 4, we can assume that $l: W \to B$ is a $4$-equivalence. Then by \cite {Kreck}, Theorem 3 the surgery obstruction $\theta (W,l)$ in $l_8(1)= l_0(1)$ is represented by 
$$
(H_4(W,M) \leftarrow KH_4(W) \to H_4(W,M'), \lambda),
$$
where $KH_4(W)$ is the kernel of $l_*:H_4(W) \to H_4(B)$, and $\lambda $ is the intersection form, which by Lefschetz duality induces a unimodular pairing $\lambda: H_4(W,M) \times H_4(W,M') \to \mathbb Z$. The maps are induced by the inclusions. In general we would have to add a quadratic refinement $\mu$ to this, but  our situation makes it true that this is determined by $\lambda$. We want to show that our conditions imply that $\Theta (W,l)$ is elementary, which by Theorem 3 of \cite {Kreck} implies that $W$ is bordant to an $h$-cobordism. Then the $h$-cobordism theorem finishes the proof. Elementary means that there is a submodule $U \subset KH_4(W)$ such that $U$ maps to direct summands in $H_4(W,M)$ and $H_4(W,M')$ of half rank and the form $\lambda$ vanishes identically on $U$.

We discuss the closed case first and assume now that $M$ and $M'$ are closed manifolds. We explain the strategy of the proof. We first note that if - as in \cite{Kreck}, Theorem 6 - $H^4(M)$ is a {\bf torsion group}, then the argument there goes through with no changes. What plays an essential role is that the map $H^4(B) \to H^4(M)$ is surjective and that $H^4(B)$ is torsion free. Both are fulfilled in our situation. Thus the three algebraic properties a) - c) on the bottom of page 749 of \cite{Kreck} are fulfilled and from that one can use the purely algebraic Proposition 13  of  \cite{Kreck} to show that the surgery obstruction is elementary. 

If $H^4(M)$ is not a torsion group we will reduce our situation to this by the following idea. We want to kill $H_4(W,M)$ and simultaneously $H_4(W,M')$ by a sequence of surgeries. The boundary operator $d: H_4(W,M) \to H_3(M)$ is surjective since the map $W \to B$ is a $3$-equivalence and this implies $H_3(W)=0$. Let $x_1,\ldots,x_r$ be in $H_4(W,M)$ mapping to a basis of $H_3(M)/_{\text{torsion}}$. The idea is to find dual (meaning that the intersection numbers fulfill $\lambda (x_i,y_i) = \delta _{i,j}$) classes $y_i$ class in $\pi_4(W)$ on which one can do surgery, which means that the self intersection number and the evaluation of $p_1(W)$ on them vanishes)  to kill them and at the same time the dual classes $x_i$ (and simultaneously for $M'$ instead of $M$). If we can find such classes, after surgeries on them we only have to kill classes in $H_4(W,M)$ and $H_4(W,M')$ that map to torsion under the boundary operator so that morally we are in the situation where $H_3(M)$ is a torsion group. 

Instead of proceeding precisely like that (which one could do) we show the existence of such classes in $\pi_4(W)$ and use the submodule $U_1$ generated by them, as a direct summand of the submodule $U$ that we have to find to show that the surgery obstruction is elementary.

Here are the details. We will construct $U$ in two steps. The normal $2$-smoothings give a surjection $H^4(B) \to H^4(M)$ and $H^4(B) \to H^4(M')$ which commutes with the isomorphism $h: H^4(M') \to H^4(M)$ (this uses the compatibility with the Pontrjagin class for the factor $BSpin$ in $B$ and that the polarization of $M'$ is the pull back of the polarization of $M$ under $(g,h)$). Dualizing $h$ we obtain an isomorphism $h^*: H_4(M) \to H_4(M')$ commuting with the map induced by the normal $2$-smoothing on  $M$ and $M'$.  That $h^*$ is an isomorphism follows since $H_4(M)$ is torsion free, which holds since by Poincar\'e duality $H_4(M) \cong H^3(M)$ and by the universal coefficient theorem the torsion of $H^3(M)$ is isomorphic to the torsion of $H_2(M)$, which is zero. Thus we have  a commutative diagram

\[
\xymatrix@=30pt{
H_4(M) \ar@{^-->}[d]^{h^*} \ar@{^{(}->}[r]^{\bar \nu _*} & H_4(B) \ar@{^-->}[d]^{=}  \\
H_4(M') \ar@{^{(}->}[r]^{(\bar \nu')_*} & H_4(B)  \\}
\]
The injectivity follows from the surjectivity in cohomology since $H_4(M)$ and $H_4(M')$ are torsion free, which holds since $Torsion \,H_4(M) \cong Torsion \, H^3(M) \cong Torsion \, H_2(M) =0$. Thus $\bar \nu_* -  \bar \nu'_*h^*: H_4(M) \to H_4(B)$ maps to zero. This implies that the image of $(i_M)_* - (i_{M'})_* h^*$ is contained in $KH_4(W)$. Here $i_M$ and $i_{M'}$ are the inclusions into $W$ from $M$ and $M'$ respectively.

$M$ resp. $M'$ to $W$. We denote the image of this map by $U_1 \subset KH_4(W)$. For later use we note that the sum of $U_1\subset H_4(W)$  with the image of $H_4(M)$ under $(i_M)_*$  is equal to the sum of $U_1$ with the image of $H_4(M')$ under $(i_{M'})_*$ and that this sum is the image of $H_4(\partial W)$ in $H_4(W)$. Since the  intersection form vanishes on classes coming from the boundary, the form $\lambda$ is zero on $U_1$. The image of $U_1$ in $H_4(W,M)$ and $H_4(W,M')$ is a direct summand. This follows from Poincar\'e Lefschetz duality, since if $x \in H_4(M)$ is a primitive element, there is $y \in H_3(M)$ with $x \cdot y = 1$. But since $H_3(W) =0$, there is a relative class $z \in H_4(W,M)$ with $d(z) = y$ and $i_*(x) \cdot z =x \cdot y =1$ implying that $i_*(x)$ is a primitive class in $H_4(W,M')$. Moreover it follows that  $U_1$ maps injectively to  $H_4(W,M)$ and to $H_4(W,M')$.  Our aim is to extend $U_1$ by a direct summand $U_2$ such that $U_1 \oplus U_2$ is $U$ with the properties we want. By construction the intersection number of elements of $U_1$ with all elements of $KH_4(W)$ is zero. Thus we only have to find $U_2$ such that the form $\lambda$ vanishes on $U_2$ and $U= U_1 \oplus U_2$ is a direct summand of half rank in $H_4(W,M)$ and $H_4(W,M')$. 

We are now able to reduce our situation to a situation which is equivalent to the argument in the proof of Theorem 6 in \cite{Kreck} and finish the proof with the arguments there. For this we consider the surjective composition $H_4(W,M) \to H_3(M) \to H_3(M)/_{Torsion(H_3(M))}$ and denote its kernel divided by the image of $U_1$ by $\widehat H_4(W,M): = \text {kernel} \,(H_4(W,M) \to H_3(M)/_{Torsion(H_3(M))})/_{\text {image}\,U_1}$. Similarly we define $\widehat H_4(W,M')$. Since $U_1$ is a  direct summand, $\widehat H_4(W,M)$ and $\widehat H_4(W,M')$ are torsion free groups. Next we recall that $H_4(M) \to H_4(B)$ is injective and its image is a direct summand and coincides with the corresponding map $H_4(M')\to H_4(B)$ under the identification of $H_4(M)$ with $H_4(M')$. Thus the image for $M'$ instead of $M$ is the same. We divide by this image and denote the quotient by $\widehat H_4(B)$, a torsion free abelian group. Similarly we divide $H_4(W)$ by the image of the sum of $H_4(M)$ and $U_1$ and denote the result by $\widehat H_4(W):= H_4(W) /_{\text {image} (H_4(M)) + U_1} = H_4(W) /_{\text {image} (H_4(M')) + U_1} = H_4(W)/_{\text {image} (H_4(\partial W))}$. The map $\widehat H_4(W) \to \widehat H_4(B)$ is surjective.  We denote its kernel by $\widehat KH_4(W)$.  Finally we denote the torsion of $H_3(M)$ by $\widehat H_3(M)$ and similarly for $H_3(M')$. Thus we have short exact sequences:
$$
0 \to \widehat H_4(W) \to \widehat H_4(W,M) \to \widehat H_3(M) \to 0
$$
and
$$
0 \to \widehat H_4(W) \to \widehat H_4(W,M') \to \widehat H_3(M') \to 0.
$$
Since the image of $H_4(M)$ and $H_4(M')$ is contained in the radical of $\lambda$, we obtain an induced form $\widehat \lambda$ on $\widehat H_4(W)$, which induces a unimodular pairing
$$
\widehat H_4(W,M) \times \widehat H_4(W,M') \to \mathbb Z
$$
Thus we obtain an element $\widehat \theta (W,l)$ in $l_8(1)$:
$$
(\widehat H_4(W,M) \leftarrow \widehat KH_4(W) \to \widehat H_4(W,M'), \hat \lambda).
$$
If $\hat \Theta (W,l)$ is elementary there is a subgroup $\widehat U_2 \subset \widehat K H_4(W) $ on which $\hat \lambda$ vanishes such that the images in $\widehat H_4(W,M)$ and $\widehat H_4(W,M')$ is a direct summand of half rank. By construction there are exact sequences
$$
0 \to U_1 \to KH_4(W) \to \widehat KH_4(W) \to 0
$$
and
$$
0 \to \widehat KH_4(W) \to \widehat H_4(W) \to \widehat H_4(B) \to 0.
$$
All groups in these exact sequences are torsion free (use the second exact sequence and the fact that $\widehat H_4(B)$  and  $\widehat H_4(W)$ are torsion free to see that $\widehat KH_4(W)$ is torison free and then the statement for the first exact squence follows). Thus there is a splitting $s$ of the first exact sequence and define $U_2 := s(\widehat U_2)$. The sum of $U_1$ and $U_2$  is the subgroup $U$ we are looking for to prove that $\Theta (W,L)$ is elementary.

Thus we are in the situation of the proof of Theorem 6 in \cite{Kreck}. The only difference is that instead of the homology groups  we have the corresponding ``roof-groups": $H_4(W)$ is replaced by $\widehat H_4(W)$, $H_3(M)$ is replaced by $\widehat H_3(M) := \text{torsion}\,\, (H_3(M))$ and $H_4(B)$ is replaced by $\widehat H_4(B)$. We also introduce $\widehat H_4(W, \partial W) := \text {kernel} (H_4(W,\partial W ) \to H_3(\partial W)/_{\text{torsion}}$ and $\widehat H_4(M) = \widehat H_4(M') := 0$. With these definitions one has exact pair and triple sequences:
$$
0=\widehat H_4(M) \oplus \widehat H_4(M') \to \widehat H_4(W) \to \widehat H_4(W,\partial W) \to \widehat H_3(M) \oplus \widehat H_3(M') \to 0
$$
$$
0 = \widehat H_4(M') \to \widehat H_4(W, M) \to \widehat H_4(W,\partial W) \to \widehat H_3(M') \to 0
$$
$$
0 = \widehat H_4(M) \to \widehat H_4(W, M') \to \widehat H_4(W,\partial W) \to \widehat H_3(M) \to 0.
$$
Then we define the corresponding roof-cohomology groups: $\widehat H^4(M) := \text{torsion}\,\, (H^4(M))$, $\widehat H^3(M) =0$ and similarly for $M'$. We define $\widehat H^4(W):= \text {kernel} (H^4(W) \to H^4(\partial W)/_{\text{torsion}})$ and $\widehat H^4(W, \partial W)/_{\text {image} H^3(\partial W)}$. We define $U^1 $ as the image of the map $H^3(M) \to H^4(W, \partial W) $ defined by the Poincar\'e dual of the map used to define $U_1$ and with this we define $\widehat H^4(W,M) := \text {kernel}(H^4(W, M) \to H^5(M)/_{\text{torsion}})/_{\text {image} (U^1)}$ and similarly for $M'$. Finally we define $\widehat H^4(B) := \text {kernel} (H^4(B) \to H^4(M))$. There are corresponding pair sequences for the roof cohomology groups as above for the roof homology groups. 

By construction Poincar\'e Lefschetz duality induces duality isomorphisms for the roof groups. There are also induced Kronecker isomorphisms. One has an induced cup product pairing and evaluation on $[W,\partial W]$ leading to  $\widehat H^4(W, \partial W) \otimes \widehat H^4(W, \partial W) \to \mathbb Z$. The conditions that $\text {sign}  (W) =0$ together with the fact that we divided by subgroups in the radical of the intersection form of $W$ implies that the signature of this form is trivial. Furthermore we consider the map $\hat l : \widehat H^4(B) \to \widehat H^4(W)$ induced by $l^*$, and the condition  $l^*(x) \cup l^*(y) =0$ for $x,y \in H^4(B) \otimes \mathbb Q$ for all $x, y \in H^4(B)$ mapping to $0$ in $H^4(\partial W) \otimes \mathbb Q$ implies $\hat l (x) \otimes \hat l(y) =0$ for all $x, y \in \widehat H^4(B) \otimes \mathbb Q$. Thus the conditions of Theorem 6 in \cite{Kreck} are fulfilled if we replace cohomology groups everywhere by the roof-groups. Our situation is different, since $B$ is not $K (\mathbb Z^m,2)$ as in Theorem 6. But -- as mentioned at the beginning of the proof -- all arguments go  through identically, where it is important that $\widehat H_4(B)$ is torsion free and the map $\widehat H_4(W) \to \widehat H_4(B)$ is surjective. We also use the identification of $H_3(M)$ with $H_3(M')$ which by Poincar\'e duality is given by the identification of $h: H^4(M)\to H^4(M')$. This induces an identification $\widehat H_3(M) \cong \widehat H_3(M')$, the group called $H$ in the proof of Theorem 6. Thus the three algebraic properties a) -- c) on the bottom of page 749 in \cite{Kreck} are fulfilled and from that one can use the purely algebraic Proposition 9  of \cite{Kreck} to show that the surgery obstruction $\widehat \theta (W,l)$ is elementary. 

This finishes the proof of Theorem 5 if $M$ and $M'$ are closed manifolds.

Now, we consider the case, where $M = N \times I$ and $\partial W = M_h$, the mapping torus of a self diffeomorphism on $N$ acting trivially on all homology groups , and $\bar \nu_{N_h}$ is a good normal structure on $N_h$. Let $(W,l)$ be -- as in the closed case -- a zero bordism, such that $l$ is a $3$-equivalence. Then we consider the map $l_*: H_4(\partial W) \to H_4(B)$. Since $\bar \nu_{N_h}$ is good, the  image is the same as the image of the precomposition with the map $H_4(M) \to H_4(\partial W)$. By the Wang sequence and the fact that $h$ acts trivially on homology groups, we have a surjection $H_4(\partial W) \to H_3(M)$. Since $H_3(N)$ is torsion free, there is a splitting and under the conditions above we can choose this splitting $s: H_3(M) \to H_4(\partial W)$ so that its image in $H_4(W)$ is contained in $KH_4(W)$. Moreover Poincar\'e duality for $\partial W$ implies that this image is a direct summand both in $H_4(W, (M \times I)_0)$ and $H_4(W, (M \times i)_1)$, where $(M \times I)_0$ and $(M \times I)_1$ are disjoint embeddings into $\partial W$ (so that the complement is the disjoint union of the mapping cylinders of $h$ and $\id$). Thus -- as in the closed case -- we define $U_1 \subset KH_4(W)$ to be the image of $s(H_3(N))$. It has the same properties as $U_1$ in the closed case. The rest of the argument is the same as in the closed case, where we replace $M$ and $M'$ by $(M \times I)_0$ and $(M \times I)_1$. Actually some points are even simpler, since in our situation $M$ and $M'$ are both homotopy equivalent to $N$, and the complements of $N$ is itself homotopy equivalent to $N$.

\end{proof}

If one wants to apply this theorem one needs to know the bordism group 
$\Omega _7(B)$. It turns out that this group is zero . 

\begin{theorem} $\Omega_7(B(n,m,w_2)) =0$.

The same holds for the corresponding topological bordism group. 
\end{theorem}

We defer the proof as well as another computation of a bordism group to section 5.

With this result one can try to classify $1$-connected closed $7$-manifolds with torsion free second homology by choosing normal $2$-smoothings in $B$ and finding conditions which allow one to modify a $B$-bordism between them by adding a closed manifold. This was done in special cases in \cite{K-S}, \cite{K-S2}. In the next section we carry this out  for $2$-connected closed $7$-manifolds reproving a Theorem by Crowley and Nordstr\"om \cite{C-N}. In the section after that we classify $7$-manifolds with same cohomology ring as $S^2 \times S^7 \sharp S^3 \times S^4$.  Other applications concern the mapping class group of simply connected $6$-manifolds $N$ with $H_2(N)$ torsion free. This will be studied in \cite{K-Y}.

\section{$2$-connected $7$-manifolds}

We begin with the definition of the quadratic refinement of the linking form 
$$b: \text{torsion}\,\, H^4(M) \otimes \text{torsion}\,\,H^4(M) \to \mathbb Q/_{\mathbb Z}
$$ and of the generalized Eells-Kuiper invariant. The definition is the following. We will show that there is a compact $spin$ $8$-manifold $W$ with signature $0$ and boundary $M$ such that the restriction $H^4(W) \to H^4(M)$ is surjective. Recall that in this situation one can define the linking form as described after Theorem 5. If $H^4(W) \to H^4(M)$ is surjective there is a class $\hat x\in H^4(W)$ restricting to the $d(M)$ structure $x$ on $M$. Then for a torsion class $y \in H^4(M)$ we choose $\hat y \in H^4(W)$ restricting to $y$ on $M$. Now we pass to rational cohomology so that $\bar p_1(W) - d(M)\hat x\in H^4(W;\mathbb Q)$ and $\hat y \in H^4(W;\mathbb Q)$ restrict to $0$ in $H^4(M;\mathbb Q)$. Thus the classes can be pulled back to $H^4(W,\partial W;\mathbb Q)$ and the following intersection number is our quadratic refinement:
$$
q_x (y):= \langle  (\bar p_1(W) - d(M) \hat x)\cup \hat y + \hat y^2, [W, \partial W] \rangle \in \mathbb Q /_{2 \mathbb Z}.
$$
By a quadratic refinement we mean that $q(y+z) - q(y) - q(z) = 2b(y,z)$. From this definition and the definition after Theorem 5 it is clear that $q_x$ is a quadratic refinement of the linking form. We will show in the proof of Theorem 1 that it is well defined. We note that our definition looks a bit different from that in \cite{C-N} but the reason is only that they divide the expression above by $2$ to obtain an invariant in $\mathbb Q/_\mathbb Z$ instead of  $\mathbb Q /_{2 \mathbb Z}$. 

The generalized Eells-Kuiper invariant is one of our $s_i$-invariants and is defined similarly to the characteristic number
$$
GEK(M,x) := \langle  (\bar p_1(W) -  d(M) \hat x)^2, [W,\partial W]\rangle  \in \mathbb Q/_{8 \cdot gcd( 28, \frac {d(M)} 2,\frac{ d(M)^2+2d(M)}8)\mathbb Z}
$$
(recall that $d(M) $ is even). Again we will show in the proof of Theorem 1, that this is well defined.

To apply Theorem 4, a special case of Theorem 5, one needs to analyze the different $d$-structures and the change of the invariants:
\begin{proposition}
The different $d(M)$ structures on $M$ are obtained from a fixed $d(M)$-structure $x$ by replacing $x$ by $x+t$ for some torsion element $t\in \text{torsion}\,\, H^4(M)$. The quadratic refinement and the Eells-Kuiper invariant for $x+t$ are:
\[q_{x+t}(y) = q_x(y) - d(M)b (t,y) \in \mathbb Q /_{2\mathbb Z} \tag {1}\]
\[GEK(M, x +t) = GEK(M,x) - 2d(M) q_x(t) +
 \frac {d(M)^2-2d(M)}2 (q_x(2t) - 2q_x(t)) \tag {2} \]
\end{proposition}

For computations it might be useful to note that  $q_x(2t) - 2q_x(t) = 2b(t,t) \in \mathbb Q/_{2\mathbb Z}$. The last equality in Proposition 7  improves the formula by Crowley and Nordstr\"om (\cite {C-N}, equation (1)). They give a transformation formula in terms of the linking form $b$ which for $d(M) =  2 \,\, mod \,\, 4$ is weaker than our formula. We will prove this proposition after the proof of Theorem 1. 

\begin{proof}[Proof of Theorem 1] For $2$-connected $7$-manifolds $M$ and $M'$ with $d(M)$-structure $x$, a polarization is  given by a resolution $\beta: \mathbb Z^m \to H^4(M)$, a surjective homomorphism. Later we will construct such a resolution with more control. Given this  we obtain the corresponding fibration $B = K(\mathbb Z^m,4) \times BSpin$ and consider a spin $B$-bordism $(W,l)$ between $(M, \bar \nu_M)$ and $(M', \bar \nu_{M'})$, which exists according to Theorem 6. By adding copies of $\pm  \mathbb HP^2$ we can assuem that  the signature of $W$ is$0$. We have to show that after adding an appropriate closed $B$-manifold to $W$ we can achieve that the conditions in Theorem 4 are fulfilled implying that $W$ is bordant to an $h$-cobordism. Thus we have to compute $\Omega _8^{Spin} (K(\mathbb Z,4))/_{\text{torsion}}$. We begin with this. 

We will show that a basis is given by the quternionic projective plane $\mathbb HP^2$ with the constant map to $K(\mathbb Z^m,4)$, which has signature $1$ and trivial $\hat A$-genus and $Bott$ (again with the constant map), a manifold with $\hat A(Bott) =1$, for example the manifold obtained from $28$ copies of the $E_8$-plumbing by gluing in $D^8$. This has signature $28\cdot 8$. We denote by $(\mathbb HP^2_i)$ the element given by $(\mathbb HP^2)$ together with the map to the $i$th factor of $K(\mathbb Z^m,4)$ given by $\bar p_1(\mathbb HP^2)$, which is a generator of $H^4(\mathbb HP^2)$), and  $(\mathbb HP^2_i)'$ the corresponding element given by $-\bar p_1(\mathbb HP^2)$.  Finally we consider for $i < j$ the manifold $(S^4 \times S^4)_{i,j}$, which stands for the product together with the two standard generators of $H^4(S^4 \times S^4)$ considered as maps to the $i$th and $j$th factor resp.. The same holds for the corresponding topological bordism group if we replace Bott by the $E_8$-manifold, an almost parallelizable topological manifold with signature $8$.

\begin{theorem} The bordism group  $\Omega_8^{Spin} (K(\mathbb Z^m,4))/_{\text{torsion}}$ is a free abelian group of rank $2m + \frac{m(m-1)}2 +2$ with basis given by:
$$
\mathbb HP^2, Bott, \mathbb HP^2_i, (\mathbb HP^2_i)', (S^4 \times S^4)_{i,j}
$$
for $1 \le i\le m$ and $1\le i<j\le m$. 

The same holds for the corresponding topological bordism group, if we replace Bott by the $E_8$-manifold.
\end{theorem}

We will prove this theorem in section 5. 

We are actually interested in the submodule of $\Omega_8^{Spin}(K(\mathbb Z^m,4))/_{\text{torsion}}$ consisting of elements with signature zero, since we only consider bordisms with signature zero between two normal $2$-smoothing. Such a bordism can be constructed from an arbitrary bordism by adding an appropriate sum of copies of $\pm \mathbb HP^2$. This subgroup has basis: $\alpha_i:= \mathbb HP^2_i - (\mathbb HP^2_i)'$, $\beta _i:= \mathbb HP^2_i - \mathbb HP^2$,  $\gamma_{i,j}:= (S^4 \times S^4)_{i,j}$ for $i < j$, and $\delta:= Bott \, - 2^5\cdot 7\mathbb HP^2$ (the signature of $Bott$ is $2^5 \cdot 7$ and  $\bar p_1(Bott) = 0$, since Bott is almost parallelizable), where in the topological case we replace by $\delta$ by $\delta ' = E_8 - 8 \mathbb HP^2$.

If $e_i$ is the basis of $H^4(K(\mathbb Z^m,4))$ given by the factors and $(Q,l)$ is a closed $B$-manifold we consider the characteristic numbers $\langle (\bar p_1(Q)- d(M)l^*( e_1))^2, [Q]\rangle $ (the special role of $e_1$ will become clear when we construct the more controlled resolution mentioned above), $\langle (\bar p_1(Q)- d(M)l^*( e_1))\cup l^*(e_i)- l^*(e_i)^2, [Q]\rangle $ for $2\le i \le m$,  and $\langle l^*(e_i)\cup l^*(e_j),[Q]\rangle $ for $2\le i\le j \le m$. We compute the values of these invariants on our generators $\alpha_1$ , $\alpha_i$ for $i >1$, $\beta_1$, $\beta_i$ for $i>1$, $\gamma_{i,j}$ for $i<j$, and $\delta$. We write the result in form of a vector, where the entries are in this order: first the evaluation on $\alpha_1$, then on $\alpha_i$ for $i>1$, then on $\beta_1$, then on $\beta_i$ for $i>1$, then on $\gamma_{i,j}$ for $i<j$, and finally on $\delta$. As usual  $(e_{i,j})$  stands for the vector which has $1$ at the place $(i,j)$ and $0$ else and  and $e_i$ is the $i$th canonical base element.

For $\langle (\bar p_1(Q)- d(M) l^*(e_1))^2, [Q]\rangle $ we obtain: 
\[ (-4d(M), (0,\ldots,0), -2d(M)+d(M)^2,(0,\ldots,0),(0,\ldots,0),8 \cdot 28), \, \tag{1} \]
for $\langle (\bar p_1(Q)- d(M)l^*(e_1)) \cup l^*(e_j) - l^*(e_j)^2, [Q]\rangle $ and $2\le j \le m$:
\[ (0, 2e_j, 0, (0,\ldots,0),(0,\ldots,0),0 ), \, \tag{2} \]
and for  $\langle l^*(e_i)\cup l^*(e_j), [Q]\rangle $ and $1\le i <j\le m$:
\[ (0, (0,\ldots,0), 0,(0,\ldots,0),( e_{i,j}),0), \, \tag{3} \]
and for  $\langle l^*(e_i)^2\, ,[Q]\rangle $ and $1\le i \le m$:
\[ (0,(0,\ldots,0),0,(e_i), (0,\ldots,0),0)\, \tag{4} \]

In the topological case we replace the first vector by 
$$
(-4d(M), (0,\ldots,0), -2d(M)+d(M)^2,(0,\ldots,0),(0,\ldots,0),8 ),
$$

With this information we show now that the generalized Eells-Kuiper invariant and the quadratic refinement of the linking form are well defined as claimed above. 
Let $W$ be a spin zero bordism of $M$ (equipped with a $d(M)$-structure) with signature $0$ inducing a surjection on the $4$-th cohomology,  then we choose, if $d(M) \ne 0$, a decomposition of $H^4(M)$ as $\mathbb Z x \oplus Q$ and of $H^4(W)$ as $\mathbb Z \hat x \oplus C$, where $\hat x$ restricts to $x$ on the boundary and $C$ maps under the restriction to $Q$. For another zero bordism $W'$ we choose a similar decomposition $H^4(W)$ as $\mathbb Z \hat x' \oplus C'$. Then, for $d(M)=0$ we choose a resolution $\mathbb Z^{m}$ to both of $H^4(W)$ and $H^4(W')$ such that the composition with the restriction to $H^4(M)$ commutes. For this choose $r$ generators $x_i$  of $H^4(W')$, and choose  elements $y_i$ in $H^4(W)$ mapping to the image of these generators in $H^4(M)$. Then choose $s$ generators $y_i$ of the kernel of $H^4(W) \to H^4(M)$. Now let $m = r+s$ and map $e_i$ for $i\le s$ to $x_i$ respectively $y_i$ and $e_{r+j}$ to $0$ respectively $z_j$.   Similarly, for  $d(M) \ne 0$  choose a resolution $\mathbb Z^{m-1}$ of both of $C$ and $C'$ as in the case $d(M) =0$  and extend it to $\mathbb Z^m$ by mapping $(1,(0,\ldots ,0))$ to   $\hat x$ respectively $\hat x'$.  This way we again obtain a resolution such that the compositions with the restrictions to $H^4(M)$ agree. Then we glue $W$ and $W'$ along $M$ to obtain a closed spin manifold $Q$ with a map to $K(\mathbb Z^m,4)$ and the difference of the characteristic numbers for $W$ and $W'$ is the corresponding characteristic number of $Q$. Looking at the first vector (1) we see that the generalized Eells-Kuiper invariant is well defined and looking at the second vector (2) we see that the quadratic refinement is well defined. In the topological case the same argument using the replacement of the first vector above we see the the reduced generalized Eells-Kuiper invariant $\bar{GEK}$ is well defined.

Now we show that the invariants in Theorem 1 determine the diffeomorphism type. For this  we choose a resolution with more control, as announced before. If $d(M) \ne 0$ we decompose $H^4(M) \cong \mathbb Zx \oplus \mathbb Z^r \oplus \,\, \oplus _{i = 1}^k \mathbb Z/_{n_i}$ and resolve it via $\mathbb Z x \oplus \mathbb Z^r \oplus \mathbb Z^k$ in the obvious way. If $d(M) =0$ we do the same omitting the first $\mathbb Z $-summand. If we have an isomorphism between $H^4(M)$ and $H^4(M')$ as in Theorem 4 we consider the induced resolution. In both cases we consider the corresponding $B$-structures $(M, \bar \nu_M)$ and $(M',\bar \nu_{M'})$, where $B = K(\mathbb Z,4) \times K(\mathbb Z^{r},4)\times K(\mathbb Z^{k},4)\times BSpin$, if $d(M) \ne 0$, and $B = K(\mathbb Z^{r},4)\times K(\mathbb Z^{k},4)\times BSpin$, if $d(M) =0$. By Theorem 6 there is a $B$-bordism $(W,l)$ with signature $0$ (after adding copies of $\pm \mathbb HP^2$)  between $(M, \bar \nu_M)$ and $(M',\bar \nu_{M'})$. We have to show that all characteristic numbers occurring in our  Theorem 4 vanish after adding appropriate closed manifolds. Then  Theorem 4 implies that we can replace $W$ by an $h$-cobordism. These are the characteristic numbers
\begin{enumerate}
\item
 $$\langle (\bar p_1(W)- d(M) l^*( x))^2, [W, \partial W]\rangle,$$
 where $x$ is the generator of the first factor $K(\mathbb Z,4)$, if $d(M) \ne 0$, and $0$ else. 
 \item 
 $$\langle (\bar p_1(W)- d(M) l^*( x))  l^*(e_i), [W, \partial W]\rangle,$$
 where again $x$ is the generator of the first factor $K(\mathbb Z,4)$, if $d(M) \ne 0$, and $0$ else, and $e_i$ is the basis of $H^4(K(\mathbb Z^m,4))$ given by the factors (these are the generators of the summand $\mathbb Z^k$ in our resolution which resolves the torsion subgroup). 
 \item
$$\langle l^*(e_i) \cup l^*(e_j), [W, \partial W]\rangle $$
 for $1 \le i \le j\le k$.
 \item 
  We actually replace $\langle (\bar p_1(W)- d(M)l^*x) l^*(e_i) , [W, \partial W]\rangle $ by 
 $$\langle (\bar p_1(W)- d(M)l^*(x)) l^*(e_i)- l^*(e_i)^2 , [W, \partial W]\rangle,$$ which  of course makes no difference for our condition, since we also assume $\langle l^*(e_i)^2, [W, \partial W]\rangle  =0$, but $\langle (\bar p_1(W)- d(M)x) l^*(e_i)- l*(e_i)^2 , [W, \partial W]\rangle $ is controlled by the quadratic refinement. 
 \end{enumerate}

The characteristic number $\langle (\bar p_1(W)- d(M) l^*( x))^2, [W, \partial W]\rangle $ is $0$ mod $\gcd(-4d(M),-2d(M)+d(M)^2,8 \cdot 28)$ since we assume that the generalized Eells-Kuiper invariants agree for $(M,x)$ and $(M',x')$. Thus, if $d(M) \ne 0$,  the first vector (1) in the computation of our invariants shows that we can make $\langle (\bar p_1(W)- d(M) l^*( x))^2, [W, \partial W]\rangle $ equal to $0$ by adding a linear combination of the $\alpha_1$ and $\beta_1$ and $\delta$.  If $d(M) = 0$,  we just add a multiple of $\delta$. 

The characteristic numbers $\langle (\bar p_1(W)- d(M)l^*(x)) l^*(e_i)- l^*(e_i)^2 , [W, \partial W]\rangle $ for $1 \le i \le k$ are $0$ mod $2$, since we control the quadratic refinement. Then the second vector (2) giving the values of these numbers on our generators implies that by adding a linear combination of the $\alpha_i$'s and $\beta_i$'s for $i>1$ we can kill these numbers.  Note that the first characteristic number is unchanged.

Finally we consider the the characteristic numbers $\langle l^*(e_i) \cup l^*(e_j), [W, \partial W]\rangle $ for $1 \le i \le j\le k$. They agree mod $\mathbb Z$ since we control the linking form. The third (3) and fourth vector (4)  giving the values of these numbers on our generators implies that by adding a linear combination of $\gamma_{i,j}$ we can kill these numbers. Note that the two classes of characteristic numbers, which we made zero before, are again not changed. 
Thus all characteristic numbers occurring in  Theorem 4 are $0$ after these changes, implying the statement. 

In the topological case the only change is that we have to replace $\delta$ by $\delta '$. 

\end{proof}

\begin{proof}[Proof of Proposition 7] 
Now we prove the formula for the change of the quadratic refinement and the generalized Eells-Kuiper invariant if we replace $x$ by $x+t$ for some torsion element in $H^4(M)$. Let $W$ be a spin zero bordism of $(M)$ with signature $0$ inducing a surjection in cohomology. We choose $\hat t \in H^4(W)$ extending  $t$.

Putting this into the formula for $q_x$ we obtain:
$$
q_{x'}(y) = q_{x +t}(y) = \langle ((\bar p_1(W) - d(M)(\hat x + \hat t)) \cup \hat y + \hat y^2), [W,\partial W]\rangle  =$$
$$
\langle (\bar p_1(W) - d(M)\hat x) \cup \hat y + \hat y^2, [W, \partial W]\rangle  - d(M)\langle \hat t \cup \hat y , [W , \partial W]\rangle  = 
$$
$$q_x(y) - d (M)\langle \hat t \cup \hat y), [W, \partial W]\rangle  \in \mathbb Q /_{2\mathbb Z}.
$$
The last expression is $d(M)b(t,y) \in \mathbb Q /_{2\mathbb Z}$. This makes sense in $ \mathbb Q /_{2\mathbb Z}$, even though $b(t,y)$ is only  well-defined in $ \mathbb Q /_{\mathbb Z}$, since $d(M)$ is even. Thus
$$
q_{x'}(x +t )(y) = q_x(y) - d(M)b (t,y) \in \mathbb Q /_{2\mathbb Z}
$$

Putting the expression into the formula for $GEK$ we obtain:
$$
GEK (M,x') = \langle ((\bar p_1(W) - d(M)\hat x) - d(M) \hat t)^2,[W,\partial W]\rangle  =  
$$
$$\langle (\bar p_1(W) - d(M)\hat x)^2,[W,\partial W]\rangle  - \langle 2d(M)(\bar p_1(W)-d(M)\hat x) \cup \hat t,[W,\partial W]\rangle  + 
$$
$$d(M) ^2\langle \hat t^2,[W, \partial W]\rangle =
$$
$$
GEK(M,x) - 2d(M) q_x(t) + \frac {d(M)^2-2d(M)}2 (q_x(2t) - 2q_x(t))
\in \mathbb Q/_{8 \cdot gcd( 28, \frac {d(M)} 2,\frac{ d(M)^2+2d(M)}8)\mathbb Z}.$$
\end{proof}

\begin{remark} We explain the relation to Theorem 1.3, second sentence, by Crowley and Nordstr\"om. They consider the set of choices for classes $x$ in our Theorem 1 and denote it by $S_{d_\pi}$. For each $x$ we defined a quadratic refinement $q_x$ which agrees with their quadratic refinement and a generalized Eells-Kuiper invariant which agrees with their Eells-Kuiper invariant. Thus one obtains an invariant of $M$ by considering this function on $S_{d_\pi}$. With this their statement is equivalent to our statement. 
\end{remark}

\section{Manifolds that have the cohomology ring of $S^2 \times S^5 \sharp S^3 \times S^4$}
We consider closed simply connected spin 
$7$-manifolds $M$ with the cohomology ring of $S^2 \times S^5 \sharp S^3 \times S^4$. Again we consider $d(M)$, the divisibility of $\bar p_1(M)$ and recall that $d(M) $ is even. We choose generators $x \in H^2(M)$ and $y \in H^4(M)$ such that $\bar p_1(M) = d(M) y$. Note that $x$ is unique up to sign and $y$ is determined by $\bar p_1(M)$, if this is non-zero. The data $x$ and $y$ determine a polarization and using the characteristic numbers occurring in Theorem 4 we define Kreck-Stolz type invariants 
$$
s_1(M,x,y) \in \mathbb Z/_{8 \cdot gcd( 28, \frac {d(M)} 2,\frac{ d(M)^2+2d(M)}8)\mathbb Z}
$$
$$
s_2(M,x,y) \in \mathbb Z/_{gcd(24,d(M))}
$$
$$
s_3(M,x,y) \in \mathbb Z/_2
$$
The definition is as follows. There is a compact oriented spin $8$-manifold $W$ with signature zero and with boundary $M$ and classes $\hat x \in H^2(W)$ restricting to $x \in H^2(M)$ and $\hat y \in H^4(W)$ restricting to $y \in H^2(Y)$. Then the relative characteristic numbers $\langle(\bar p_1(W) - d(M)\hat y) ^2, [W,\partial W]\rangle  $, $\langle (\hat x^2 )^2, [W,\partial W]\rangle $, and $\langle \hat x^2 (\bar p_1(W) - d(M)\hat y),[W,\partial W]\rangle $ are integers, since the restriction to the boundary of at least one of the factors is zero, allowing us to pull them back to $H^4(W,M)$ and taking the cup product with the absolute class. From these we define the invariants $s_i$ by the formulae
$$
s_1(M,x,y) := [\langle (\bar p_1(W) - d(M)\hat y) ^2,[W,\partial W]\rangle ] \in  \mathbb Z/_{8 \cdot \gcd( 28, \frac {d(M)} 2,\frac{ d(M)^2+2d(M)}8)\mathbb Z}$$
$$
s_2(M,x,y) :=  [ \langle (\hat x^2)^2  - \hat x^2 ( \bar p_1(W) - d(M) \hat y), [W,\partial W]\rangle ]
\in \mathbb Z/_{\gcd (24,d(M))}
$$
and 
$$
s_3(M,x,y) := [\langle \hat x^2( \bar p_1(W) - d(M) \hat y), [W,\partial W]\rangle ]\in \mathbb Z/_2.
$$
We will show that these invariants are well defined. We note that if we change the sign of $x$, the invariants are unchanged. If $\bar p_1(M) \ne 0$, then $y$ is determined by $\bar p_1(M)$, and if $d(M) = 0$, the invariants are independent on the sign of $y$. Thus in all cases we obtain well defined invariants $s_i(M) := s_i(M,x,y)$. 

Now we are ready to prove Theorem 2 (including the proof that the invariants above are well defined).

\begin{proof} (Theorem 2)
The generators $x_M$ and $y_M$ can be interpreted as maps to $K(\mathbb Z,2) \times K(\mathbb Z,4)$ giving a polarization as needed for  Theorem 4. 

We first show that the invariants $s_i$ are well defined. For this we have to compute $\Omega_8 ^{Spin}(K(\mathbb Z,2) \times K(\mathbb Z,4))/_{\text{torsion}}$, which we will do in section 5.

\begin{theorem} The subgroup of $\Omega_8 ^{Spin}(K(\mathbb Z,2) \times K(\mathbb Z,4))/_{\text{torsion}}$ of elements with signature $0$ has  basis:
$$ 
V_1 := ((Bott,0,0) - 2^5  7 \mathbb HP^2, \,\, V_2:= (V(2),x,0) - 2 \mathbb HP^2,\,\, V_3 :=((S^2)^4, \Delta,0), 
$$
$$V_4 := (\mathbb HP^2, 0,y) -  \mathbb HP^2, \,\,V_5 := (S^4 \times S^4, 0,\Delta), V_6 = \frac 1 2 (S^2 \times S^2 \times S^4, \Delta ,y)
$$
\end{theorem}

Here $Bott$ is the Bott manifold obtained from the boundary connected sum of $28$ copies of the $E_8$ plumbing in dimension $8$ with an $8$-disc attached to the boundary (there are two diffeomorphisms on $S^7$ one can use to attach the $8$-disc, the spin bordism class is not affected by this choice). $V(2) \subset \mathbb CP^5$ is the degree $2$ hypersurface and $x$, $y$ are generators of the second resp. 4-th cohomology groups,  $\Delta$ stands in all cases for the diagonal class. For $V_6$ we will prove that the bordism class $(S^2 \times S^2 \times S^4, \Delta ,y)$ is divisible by $2$. 

We have to determine the value of the characteristic numbers occurring in  Theorem 2. These are (recall that $d=2s$)
$$
 (\bar p_1(V_i) - 2s\hat y) ^2, \,\, ((\hat x^2)^2  - \hat x^2 \bar p_1(V_i) + 2s\hat x^2 \hat y), \,\,
\hat x^2( \bar p_1(V_i) - 2s \hat y).
$$

The values of the first characteristic number for the $(V_i,x,y)$ are:
\[(-2^57, 0,0, 4s^2 -4s,8s^2,0), \, \tag{5} \]
for the second characteristic number
\[(0, 0,24, 0,0,0-2s), \, \tag{6} \]
and for the third
$$
(0, 2,0, 0,0-2s),
$$
After changing the base given by the $(V_i,x,y)$ by replacing $(V_6,x,y)$ by the sum of $(V_6,x,y)$ with $s$ copies of $(V_2,x,y)$, we replace the last vector by
\[(0,2,0,0,0,0). \, \tag{7} \]
This implies that the invariants $s_i$ in Theorem 2 are well defined. 

In addition, if these invariants agree for $(M,x_M,y_M)$ and $(M',x_{M'},y_{M'})$, there is a bordism $(W,\hat x, \hat y)$ with signature of $W$ zero between $(M,x_M,y_M)$ and $(M',x_{M'},y_{M'})$, such that the characteristic numbers 
$$
 (\bar p_1(W) - 2s\hat y) ^2, \,\, ((\hat x^2)^2  - \hat x^2 \bar p_1(W) + 2s\hat x^2 \hat y), \,\,
\hat x^2( \bar p_1(W) - 2s \hat y)
$$
take values in the lattice given by their values on the closed manifolds $(V_i,x,y)$. Thus after adding the disjoint union of an appropriate linear combination of these closed manifolds to $W$, we shall assume that the characteristic numbers vanish for $W$. But these characteristic numbers generate all the characteristic numbers occurring in Theorem 4 finishing the proof of Theorem 2.
\end{proof}

Now we prove Corollary 3:

\begin{proof}[Proof of Corollay 3] 
We apply Theorem 2 to classify the manifolds $M_{s,k} \sharp \Sigma$. These manifolds are the boundary of the manifold $Q$ obtained as boundary connected sum of $W_s$, the disc bundle associated to the sphere $S_s$ with $P_k$, the disc bundle associated to $T_k$, and $r$ copies of the plumbing of the $E_8$-lattice in dimension $8$, whose boundary is the homotopy sphere $\Sigma _r$. For $r=1\ldots27 \,\, mod \,\, 28$ these are the $27$ exotic $7$-spheres (for  $r= 0 \,\, mod \,\, 28$ one obtains the standard sphere). The manifold $Q$  is a spin manifold and the generators $x$ and $y$ on $M_{s,k} \sharp \Sigma$ extend to classes $\hat x \in H^2(W_s) \cong H^2(Q)$ and $\hat y \in H^4(T_k) \cong H^4(Q)$. The signature of $W_s$ is zero, of $P_k$ is $1$ and of the $r$ copies of the plumbing construction is $8r$. Thus we have to add $8r+1$ copies of $- \mathbb HP^2$, to obtain a manifold with signature $0$. The characteristic numbers are
$$
\langle (\hat x ^2) ^2, [Q]\rangle  = 1$$
$$
\langle  \hat x ^2 \cup  (\bar p_1(Q) - 2s\hat y) , [Q]\rangle = \langle  \hat x^2 \cup \bar p_1(P_k)[P_k]\rangle  = 3 + 2k
 $$
 $$ \langle (\bar p_1(Q) - 2s\hat y)^2 ,[Q]\rangle = \langle \bar p_1(P_k)^2),[P_k]\rangle - (8r+1) \langle \bar p_1(\mathbb HP^2)^2 ,[\mathbb HP^2] \rangle =  8(1-r) + 6k + 4k^2 
$$
Thus 
$$
s_1( M_{s,k} \sharp \Sigma_r, x,y) = 8(1-r)  + 6k + 4k^2 \,\, mod 8\,\,  \gcd (28,\frac{s(s-1)}2,s^2)
$$
$$
s_2( M_{s,k} \sharp \Sigma_r, x,y) = 1-2k \,\, \mod \,\, \gcd (24,2s)
$$
$$
s_3( M_{s,k} \sharp \Sigma_r, x,y) = 1 \,\, \mod\,\, 2
$$
This finishes the proof of Corollary 3.
\end{proof}

\section{Computations of some bordism groups}

Proof of Theorem 6: $\Omega_7(B(n,m,w_2)) =0$. The same holds for the corresponding topological bordism group. 

\begin{proof}[Proof of Theorem 6]  The proof in the topological and smooth category is the same, since the Topspin bordism groups are isomorphic to the spin bordism groups under the forgetful map from spin manifolds to Topspin manifolds in degree $ \le 7$ except in degree $4$, where both groups are $\mathbb Z$ but the forgetful map is multiplication by $2$ (\cite{K-S2}, Lemma 6.4). Thus we restrict ourselves to the smooth case.

We use the Atiyah-Hirzebruch spectral sequence or, if $w_2 \ne 0$, the James spectral sequence \cite{T}. The entries in the $E_2$-tableau are $H_p(K(\mathbb Z^n,2) \times K(\mathbb Z^m,4); \Omega _q^{Spin})$. The non-trivial reduced homology groups in degree $\le 7$ of $K(\mathbb Z,4)$ are $\mathbb Z$ in degree $4$ and $\mathbb Z/2$ in degree $6$ \cite {E-M}. The Steenrod square $H^4(K(\mathbb Z,4); \mathbb Z/2) \to H^6(K(\mathbb Z,4);\mathbb Z/2)$ is non-trivial \cite {E-M}. Thus the only non-trivial entry relevant for degree $7$ in the $E^2$-tableau is $H_6(B; \Omega_1^{Spin}) = H_6(B; \mathbb Z/2)$. The $d_2$-differential $H_6(B;  \Omega_1^{Spin}) \to H_4(B; \Omega _2^{Spin})$ is dual to the map $ x \mapsto Sq^2(x) + w_2 \cup x$  and similar for the differential ending in $H_6(B;  \Omega_1^{Spin}) $ \cite{T}. 

Now we abbreviate $K(\mathbb Z^n,2)$ by $K$ and $K(\mathbb Z^m,4)$ by $K'$. For the next argument we note, that since $w_2 \in H^2(K)$ there is an inclusion $B(n,0,w_2) \to B(n,m,w_2)$ and a projection $B(n,m,w_2) \to B(n,0,w_2)$. This implies that we can splitt off the spectral sequence for computing $\Omega_7(B(n,0,w_2))$ from that computing $\Omega_7(B(n,m,w_2))$. The $E^3$-tableau for the James spectral sequence for computing $\Omega_7(B(n,0,w_2))$ is zero \cite{K-S}. The spitting of the spectral sequence implies that also the $E^3$ tableau for  computing $\Omega_7(B(n,m,w_2))$ is zero, if the following sequence with $\mathbb Z/2$-coefficients is exact:
$$
H^0(K) \otimes H^4(K') \stackrel {d}{\to} H^0(K) \otimes H^6(K') \oplus H^2(K) \otimes H^4(K') \stackrel{d'}{\to} 
$$
$$
H^0(K) \otimes H^8(K') \oplus H^2(K) \otimes H^6(K') \oplus H^4(K) \otimes H^4(K')
$$
where $d$ and $d'$ are dual to the $d_2$ -differentials in the $E^2$ term. We equip $H^2(K)$ with the basis $x_i$ coming from the factors and  $H^4(K')$ with the basis $y_j$ coming from the factors. Then $H^6(K')$ is equipped with the basis $Sq^2(y_j)$. Then  the differential $d$ is given by mapping
$$
1 \otimes y_j \mapsto 1 \otimes Sq^2 y_j + w_2 \otimes y_j.
$$
The differential $d'$ is given by (use that $Sq^2 Sq^2 (y_j) =0$)
$$
1 \otimes Sq^2(y_j) \mapsto w_2 \otimes Sq^2(y_j)$$
and
$$
x_i \otimes y_j \mapsto (x_i^2 + w_2 x_i)\otimes y_j +  x_i Sq^2(y_j).
$$
Now it is an exercise in linear algebra to show that this sequence is exact implying that the $E^3$-term is $0$.

\end{proof}

Now we prove Theorem 8: The bordism group  $\Omega_8^{Spin} (K(\mathbb Z^m,4))/_{\text{torsion}}$ is a free abelian group of rank $2m + \frac{m(m-1)}2 +2$ with basis given by:
$$
\mathbb HP^2, Bott, \mathbb HP^2_i, (\mathbb HP^2_i)', (S^4 \times S^4)_{i,j}
$$
for $1 \le i\le m$ and $1\le i<j\le m$. 

The same holds for the corresponding topological bordism group, if we replace Bott by the $E_8$-manifold.

\begin{proof}[Proof of Theorem 8] The proof in the topological and smooth case is the same, namely we 
will give invariants which detect the bordism classes in the bordism groups mod torsion. We abbreviate $K(\mathbb Z^m,4)$ by $K$.  Since the rational Atiyah Hirzebruch spectral sequence collapses there is an isomorphism $\Omega_n^{Spin}(K)\otimes \mathbb Q  \to \sum _{p+q = n} H_p(K; \Omega_q^{Spin}\otimes \mathbb Q)$. Thus the reduced bordism group $\Omega_8^{Spin}(K)$ mod torsion has the same rank as $H_4(K(\mathbb Z^m,4)) \oplus H_8(K(\mathbb Z^m,4))$, which is $2n + n(n-1)/2$. We will construct bordism invariants and show that they induce a surjection 
$$
\tilde \Omega_8^{Spin}(K)/_{\text{torsion}} \to \mathbb Z^{2n + n(n-1)/2}
$$
For this  let $e_i  \in H^4(K(\mathbb Z^m,4))$ be the basis which is dual to the basis given by the embeddings of the $4$-spheres generating $\pi_4(K(\mathbb Z,4))$ into the different factors. Then we consider the following invariants on a closed spin-manifold $Q$ together with a map $g: Q  \to K(\mathbb Z^m,4)$
$$
\langle g^*(e_i)^2, [Q]\rangle \in \mathbb Z
$$
$$
\frac{\langle \bar p_1(Q)\cup g^*(e_i) - g^*(e_i)^2, [Q]\rangle }2\in \mathbb Z
$$
and for $i<j$
$$
\langle g^*(e_i) \cup g^*(e_j), [Q]\rangle \in \mathbb Z.
$$
The second expression is integral since $\bar p_1(M) \,\, \mod \,\, 2$ is equal to the Wu-class (\cite{K-S2}, Lemma 6.5). There it is also proved that this holds in the topological category, so that the invariant is also defined for the corresponding topological $TOPSpin$ bordism group.

Now it is an easy exercise to show that the invariants restricted to the subgroup generated by the  elements $\mathbb HP^2_i, (\mathbb HP^2_i)', S^4 \times S^4_{i,j}
$
 yield an isomorphism to $\mathbb Z^{2n + n(n-1)/2}$. Thus this group is isomorphic to $\tilde \Omega^{Spin} (K(\mathbb Z^m,4)/_{\text{torsion}}$ and the same holds in the topological category. Passing to the unreduced bordism groups completes the statement, if we use that $\Omega_8^{Spin}$ has basis $\mathbb HP^2$ and $Bott$ and in the topological category the bordism group mod torsion  has basis $\mathbb HP^2$ and $E_8$ (see \cite {K-S2} , section 6).

\end{proof}

Now we prove Theorem 10: The subgroup of $\Omega_8 ^{Spin}(K(\mathbb Z,2) \times K(\mathbb Z,4))/_{\text{torsion}}$ of elements with signature $0$ has  basis:
$$ 
V_1 := ((Bott,0,0) - 2^5  7 \mathbb HP^2, \,\, V_2:= (V(2),x,0) - 2 \mathbb HP^2,\,\, V_3 :=((S^2)^4, \Delta,0), 
$$
$$V_4 := (\mathbb HP^2, 0,y) -  \mathbb HP^2, \,\,V_5 := (S^4 \times S^4, 0,\Delta), V_6 = \frac 1 2 (S^2 \times S^2 \times S^4, \Delta ,y)
$$

\begin{proof}[Proof of Theorem 10]  Using an inclusion of the two factors we see that the bordism group contains $\Omega_8^{Spin} (K(\mathbb Z,2))$ and $\Omega_8^{Spin} (K(\mathbb Z,4))$ as direct summands. A basis of the $\Omega_8^{Spin} (K(\mathbb Z,2))$ was computed in \cite{K-S} showing that $V_1$, $V_2$ and $V_3$ generate the subgroup given by this factor and above we computed the basis for factor corresponding to $\Omega_8^{Spin} (K(\mathbb Z,4))$ showing that $V_4$ and $V_5$ (with the base change mentioned at the end of the previous proof) generate this summand.

If we look at the Atiyah-Hirzebruch spectral sequence computing $\Omega_8^{Spin} (K(\mathbb Z,2) \times K(\mathbb Z,4))$ there is a corresponding splitting. Thus the $8$-line splits as a direct summand of the two $8$-lines for $\Omega_8^{Spin} (K(\mathbb Z,2))$ and $\Omega_8^{Spin} (K(\mathbb Z,4))$  plus the following entries: $H_2(K(\mathbb Z,2));\mathbb Z/2) \otimes H_4(K(\mathbb Z,4);\mathbb Z/2)$ sitting in $E^2_{6,2}$ and  $H_2(K(\mathbb Z,2)) \otimes H_6(K(\mathbb Z,4)) \oplus H_4(K(\mathbb Z,2)) \otimes H_4(K(\mathbb Z,4))$ sitting in $E^2_{8,0}$. The term in $E^2_{6,2}$ is killed by an incoming $d_2$ differential, since $Sq^2 : H^2(K(\mathbb Z,2);\mathbb Z/2) \otimes H^4(K(\mathbb Z,4);\mathbb Z/2) \to H^4(K(\mathbb Z,2);\mathbb Z/2) \otimes H^4(K(\mathbb Z,4);\mathbb Z/2) \oplus H^2(K(\mathbb Z,2);\mathbb Z/2) \otimes H^6(K(\mathbb Z,4)$ is injective. The same argument shows that  the $d_2$ differential starting in $H_4(K(\mathbb Z,2)) \otimes H_4(K(\mathbb Z,4)) \oplus H_2(K(\mathbb Z,2)) \otimes H_6(K(\mathbb Z,4)) = \mathbb Z \oplus \mathbb Z/2$ ending in $H_2(K(\mathbb Z,2); \mathbb Z/2) \otimes H_4(K(\mathbb Z,4); \mathbb Z/2) = \mathbb Z/2$ maps $(a,[b])$ to $[a+b]$. Thus there is an element  $(M,g)$ with $M$ a spin manifold and $g_*([M]) = (1,[1])$. On the other hand the element $(S^2 \times S^2 \times S^4, \Delta, y)$ maps under the Hurewicz homomorphism given by the image of the  fundamental class to twice the generator of $H_4(K(\mathbb Z,2)) \otimes H_4(K(\mathbb Z,4)) = \mathbb Z$. This implies that the bordism class of $(S^2 \times S^2 \times S^4, \Delta, y)$ is divisible by $2$ and together with the $V_i$ for $1 \le i \le 5$ is a basis of $\Omega_8^{Spin} (K(\mathbb Z,2) \times K(\mathbb Z,4))/_{\text{torsion}}$. 

\end{proof}

\noindent
{\small Mathematisches Institut, Universit\"at Bonn \\
and Mathematisches Institut, Universit\"at Frankfurt\\
kreck@math.uni-bonn.de}


\begin{thebibliography}{99}    
\scriptsize 

\vspace{-7pt}
\bibitem{Ce} J.Cerf, La stratification naturelle des espaces de fonctions différentiables réelles et le théorème de la pseudo-isotopie, Inst. Hautes Études Sci. Publ. Math. No. 39 (1970) 5?173.

\vspace{-7pt}
\bibitem{C} D. Crowley, The classification of highly connected manifolds in dimensions 7 and 15. Thesis (Ph.D.) Indiana University. 2002.


 \vspace{-7pt}

\bibitem{C-N}
D. Crowley and J. Nordstr\"om,
The classification of 2-connected 7-manifolds. arXiv:1406.2226 (2014)


\vspace{-7pt}

\bibitem{E-M} S. Eilenberg and S. MacLane,  On the Groups $H(\Pi, n)$, II: Methods of Computation
Annals of Mathematics
Second Series, Vol. 60, pp. 49-139

\vspace{-7pt}

\bibitem{Kreck} M. Kreck,
Surgery and duality.
Ann. of Math. (2) 149 (1999), 707 - 754.

\vspace{-7pt}

\bibitem{K-S}  M. Kreck and S. Stolz,
A diffeomorphism classification of 7-dimensional homogeneous Einstein manifolds with $SU(3)×\times SU(2)\times U(1)$-symmetry. Ann. of Math. (2) 127 (1988), 373-388

\vspace{-7pt}

\bibitem{K-S2}  M. Kreck and S. Stolz,   
Some nondiffeomorphic homeomorphic homogeneous 7-manifolds with positive sectional curvature. J. Differential Geom. 33 (1991), 465 - 486
\vspace{-7pt}

\bibitem{K-Y}  M. Kreck and Su Yang, On the mapping class group of certain $6$-manifolds. In preparation
\vspace{-7pt}
   
\bibitem{T}  P. Teichner,  
On the signature of four-manifolds with universal covering spin. Math. Ann. 295 (1993), 745 - 759
\end{thebibliography}
\end {document}